\documentclass[11pt]{article}
\usepackage{amsfonts}

\begin{document}

\def\det{\mbox{det}}
\def\Spin{\mbox{Spin}}
\def\Tr{\mbox{Tr \-}}
\def\inj{\mbox{inj \-}}
\def\R{{\mathbb R}}
\def\C{{\mathbb C}}
\def\H{{\mathbb H}}
\def\Ca{{\mathbb Ca}}
\def\Z{{\mathbb Z}}
\def\N{{\mathbb N}}
\def\Q{{\mathbb Q}}
\def\Ad{\mbox{Ad \-}}
\def\k{{\bf k}}
\def\l{{\bf l}}
\def\diag{\mbox{diag \-}}
\def\a{{\alpha}}

\begin{center}
{\Large On new explicit Riemannian SU(2(n+1))-holonomy metrics.
\footnote{This work was supported by the RFBR (grant 09-01-00142-a),
by the program "Leading Scientific Schools" (project N.
NSh-7256.2010.1), by Federal Target Grant "Scientific and
educational personnel of innovation Russia" for 2009-2013
(government contract No. 02.740.11.0429) and the Joint Project of
the Siberian and Ural Divisions of the Russian Academy of Sciences
(Grant No. 46).} }
\end{center}

\begin{center}
E.~G.~Malkovich
\end{center}

\begin{abstract}

We construct in an explicit algebraic form a family of complete
noncompact Ricci-flat metrics which generalize Calabi metrics in
real dimension $4(n+1)$ and with holonomy $SU(2(n+1))$.


\end{abstract}

\sloppy

\section[]{Introduction.}

This article is concerned with an exploration of the Ricci-flat
Riemannian metrics with exceptional holonomies and naturally
continues a number of works \cite{Baz1, Baz-Malk, Baz2, Baz-Malk2}.
In \cite{Baz-Malk2} we were studying Riemannian metrics with
$Spin(7)$-holonomy on the cones over $3$-Sasakian $7$-manifolds and
were able to find in an explicit form a continuous family of
complete noncompact $8$-dimensional metrics $\bar{g}_{\alpha}$
depending on real parameter $0\leq \alpha \leq1$. Metric $\bar{g}_0$
coincides with Calabi $SU(4)$-holonomy metric; metric $\bar{g}_1$
coincides with hyperk\"{a}hler Calabi metric with holonomy
$Sp(2)\subset SU(4)$. Every found metric
$\overline{g}_\alpha,~0<\alpha<1$ is $SU(4)$-holonomy metric and
automatically Ricci-flat. Thus, Calabi metrics are "connected" by
the obtained one-dimensional family.

However Calabi metrics (firstly appeared in \cite{Calabi}) are
both correctly defined not only for dimension $8$ but for any
dimension divisible by $4$. And a question of a generalization of
the family constructed in \cite{Baz-Malk2} for higher dimensions
is very natural. In this paper that question is positively
resolved: for any real dimension $4(n+1)$ we construct in an
explicit form the continuous family of metrics $\bar{G}_{\alpha}$
"connecting" Calabi metrics.

\vskip0.2cm

{\bf Theorem.} {\it The following family consists of complete
Ricci-flat $4(n+1)$-dimensional Riemannian metrics:
$$
\begin{array}{c}
\bar{G}_\alpha= \frac{r^4 (r^4-\alpha^4)^n
}{(r^4-\alpha^4)^{n+1}-(1-\alpha^4)^{n+1}} dr^2 +
\frac{(r^4-\alpha^4)^{n+1}-(1-\alpha^4)^{n+1}}{r^2
(r^4-\alpha^4)^n } \eta_1^2 + r^2 (\eta_2^2+ \eta_3^2)
\\ \\
+(r^2+\alpha^2) \sum_{\beta=1}^n(\eta_{4\beta}^2+\eta_{5\beta}^2)
+(r^2-\alpha^2) \sum_{\beta=1}^n(\eta_{6\beta}^2+\eta_{7\beta}^2),
\end{array}
$$
where $0 \leq \alpha \leq 1$, $r \geq 1$. Metrics $\bar{G}_0$ and
$\bar{G}_1$ have holonomies $SU(2(n+1))$ and $Sp(n+1)$ accordingly
and coincide with high-dimensional Calabi metrics found in
\cite{Calabi}. Metrics $\bar{G}_{\alpha}$ for $0<\alpha<1$ have
holonomy $SU(2(n+1))$ and for $n=1$ coincide with the family
constructed in \cite{Baz-Malk2}. For $0<\alpha<1$ metrics
$\bar{G}_\alpha$ are defined on the $(n+1)$th tensor power of the
complex line bundle over the space of complex flags in $\C^{2n+1}$
and metric $\bar{G}_1$ is defined on $T^*\mathbb{C}P^{n+1}$. }

In the next section we will explain in detail construction of the
metrics $\bar{G}_\alpha$ and will prove the above theorem.

\section[]{The Proof.}

In paper \cite{Baz-Malk2} we explored the existence of the
$8$-dimensional metrics with holonomy $Spin(7)$ of the following
form
$$
dt^2+ A_1(t)^2 \eta_1^2 + A_2(t)^2 \eta_2^2 + A_3(t)^2 \eta_3^2+
B(t)^2 (\eta_4^2+\eta_5^2) + C(t)^2 (\eta_6^2+\eta_7^2)
$$
on the cone over $7$-dimensional $3$-Sasakian manifold $M$, whose
$4$-dimensional quaternionic K\"{a}hler orbifold ${\cal O}$
possesses a K\"{a}hler structure. The form $dt$ corresponds to the
generator of the cone, the forms $\eta_i,~i=1,2,3$ are the
characteristic forms of the $3$-Sasakian manifold and the forms
$\eta_i,~i=4,\ldots,7$ are $1$-forms on the orbifold. The condition
for holonomy to be contained in $Spin(7)$ is equivalent to some
system of ODE's on the functions $A_i,B,C$. This system was explored
carefully and we were able to find the following family of
solutions:
$$
\begin{array}{c}
\bar{g}_\alpha= \frac{r^4 (r^2-\alpha^2) (r^2 +\alpha^2)}{r^8-2
\alpha^4 (r^4-1) -1} dr^2 + \frac{r^8-2 \alpha^4 (r^4-1) -1}{r^2
(r^2-\alpha^2) (r^2 +\alpha^2)} \eta_1^2 + r^2 (\eta_2^2+
\eta_3^2)
\\ \\
+ (r^2+\alpha^2) (\eta_4^2+\eta_5^2) + (r^2-\alpha^2)
(\eta_6^2+\eta_7^2),
\end{array}\eqno{(*)}
$$
where $0 \leq \alpha \leq 1$ and $r\geq 1$. Metrics $(*)$ are
defined on a smooth manifold if and only if $M$ is the Aloff-Wallach
space $N_{1,1}=SU(3)/S^1$.

Calabi constructed his metrics $\bar{G}_0$ and $\bar{G}_1$ on the
complex bundles over the K\"{a}hler-Einstein manifolds
\cite{Calabi}. In particular, metrics with holonomy $SU(n)$ were
constructed on the line bundles over the compact K\"{a}hler-Einstein
manifolds and the hyperk\"{a}hler metrics were constructed on the
$T^*\mathbb{C}P^n$. Nevertheless in \cite{Calabi} this metrics were
not written out in the explicit form.

An expression for the metric $\bar{G}_0$ was found in
\cite{Page-Pope}:
$$
[1-(\frac{1}{\rho})^{2m+2}]^{-1}d\rho^2
+[1-(\frac{1}{\rho})^{2m+2}]\rho^2(d\tau-2A)^2
+\rho^2ds^2,\eqno{(1)}
$$
where $ds^2$ is a metric on a $3$-dimensional Hodge and
K\"{a}hler-Einstein manifold $F$, $dA$ --- K\"{a}hler form on $F$.
For $m=3$ and $F=SU(3)/T^2$ metrics $(1)$ and $\bar{g}_0$ coincide.
Metric $(1)$ is defined on the $(m+1)$th power of the canonical line
bundle over $F$.

In \cite{CGLP} it was attempted to explore the metrics of
cohomogeneity one on the $T^*\C P^{(n+1)}$. It's not difficult to
understand the spherical subbundle in $T^*\C P^{(n+1)}$ fibres over
the space $SU(n+2)/(U(n)\times U(1))$. On the Lie algebra $su(n+1)$
one can choose the cobasis of the left-invariant $1$-forms $L_A^B$
such that its exterior algebra satisfies $dL_{A}^{B}=iL_A^C \wedge
L_C^B$. Index $A$ takes values in $(1,2,\beta)$ and index $\beta$
takes values from $1$ to $n$, further $\beta$ is never fixed and no
ambiguity appears. Obviously, $u(n) \oplus u(1)$ is a Lie subalgebra
in $su(n+2)$ and is not an exterior subalgebra. Forms
$L_1^{\beta}=\sigma_{\beta}$, $L_2^{\beta}=\Sigma_{\beta}$ and
$L_1^2=\nu$ constitute a basis on the quotient $su(n+2)/(u(n)\oplus
u(1))$. Then one can define the real forms: $\sigma_{1\beta}+
i\sigma_{2\beta}=\sigma_\beta$, etc. Form $\lambda =L^1_1-L^2_2$ is
real by definition. In \cite{CGLP} the metrics of the following form
$$
dt^2+ a(t)^2|\sigma_\beta|^2 +b(t)^2|\Sigma_\beta|^2+
c(t)^2|\nu|^2+f(t)^2\lambda^2, \eqno{(2)}
$$
were considered, where the summation over the $\beta$ from $1$ to
$n$ is omitted. The expression for hyperk\"{a}hler metric
$\bar{G}_1$ on the $T^* \C P^{n+1}$ was found:
$$
\frac{dr^2}{1-r^{-4}}+\frac{1-r^{-4}}{4}r^2\lambda^2+
r^2(\nu_1^2+\nu_2^2)+\frac{r^2+1}{2}(\Sigma_{1\beta}^2+\Sigma_{2\beta}^2)
+ \frac{r^2-1}{2}(\sigma_{1\beta}^2+\sigma_{2\beta}^2).\eqno{(3)}
$$
In the case $n=1$ to make the notations of papers \cite{Baz-Malk2}
and \cite{CGLP} agree one should put
$$
\lambda=2\eta_1,~\nu_1=\eta_3,~\nu_2=\eta_2,~
\Sigma_1=\sqrt{2}\eta_4,~ \Sigma_2=\sqrt{2}\eta_5,~
\sigma_1=\sqrt{2}\eta_6,~ \sigma_2=\sqrt{2}\eta_7.
$$

In \cite{CGLP} also the Ricci tensor were written out. Obviously
this tensor has five components: $Ric=R_0dt^2 +R_a|\sigma_\beta|^2
+R_b|\Sigma_\beta|^2+ R_c|\nu|^2 +R_f\lambda^2$ and depends on four
functions. We don't write it out here.

Notice that for any dimension $n$ coefficients of metric $(3)$ have
the same form, and coefficients of $(1)$ depend on $n$ explicitly.
Therefore, the requested family of metrics should depend on $n$
explicitly and for $\alpha=1$ should coincide with $(3)$. We will
look for metrics of the following form:
$$
\frac{dr^2}{W^2}+\frac{W^2r^2}{4}\lambda^2+ r^2(\nu_1^2+\nu_2^2)+
\frac{(r^2-\alpha^2)}{2}(\sigma_{1\beta}^2+\sigma_{2\beta}^2)+
\frac{(r^2+\alpha^2)}{2}(\Sigma_{1\beta}^2+\Sigma_{2\beta}^2),
$$
where $W=W(r,n,\alpha)$ is an unknown function. If one will put
appropriate functions into the expressions for the Ricci tensor then
the components $(R_a,~R_b,~R_c)$ will be
$$
\begin{array}{ll}
R_a=-\frac{2Q}{(r^2-\alpha^2)^2(r^2+\alpha^2)}\\
R_b=\frac{2Q}{(r^2+\alpha^2)^2(r^2-\alpha^2)}\\
R_c=-\frac{2Q}{r^2(r^4-\alpha^4)},
\end{array}
$$
where
$Q=\frac{dW}{dr}(r^5-r\alpha^4)+4W^2\alpha^4+4(n+1)(r^4-\alpha^4-r^4W^2)$
and $\frac{dr}{dt}=W$. This equation can be integrated without
difficulties:
$$
W^2=\frac{(r^4-\alpha^4)^{n+1}+C}{r^4(r^4-\alpha^4)^n},
$$
where $C$ is an integration constant. By a shift along $r$ this
constant can be fixed. One should put $C=-(1-\alpha^4)^{n+1}$ then
$r\geq 1$ and $W(1)=0$. The components $R_0$ and $R_f$ are the
second order ODE's and automatically vanish.

Consider the following $2$-form
$$
\Omega=rdr\wedge \lambda+ 2r^2\nu_1\wedge \nu_2
-(r^2+\alpha^2)\Sigma_{1\beta}\wedge \Sigma_{2\beta}
+(r^2-\alpha^2)\sigma_{1\beta}\wedge \sigma_{2\beta}.
$$
Using the exterior algebra's relations from \cite{CGLP} one can
easily verify that this form is closed and up to multiplying by
$\frac{1}{2}$ is the K\"{a}hler form of metric $(2)$. The vanishing
of the Ricci tensor and the closeness of the form $\Omega$ give us
the statement about holonomy.

Next we will explore the topology of the spaces where the founded
metrics are defined. Here we generalize construction for higher
dimensions used in \cite{Baz-Malk2}. Consider the complex space
$\C^{n+2}$ and the diagonal action of a circle $S^1$ on it. This
action defines an equivalence class. We will designate such a class
by square brackets. For example, $[u]$, $[V]$ where $u,~V$ are
vectors or subspaces.

Consider the space $\tilde{H}=\{ (u_1,u_2,V)|~|u_1|=1, u_1\perp
_{\C} u_2 \perp _{\C} V \} \subset S^{2n+3}\times \C^{n+2} \times
G_n(\C ^{n+2})$. Consider also the projection $\tilde{\pi} :
(u_1,u_2,V) \rightarrow (u_1,V)$ from $\tilde{H}$ to the space
$\tilde{F}= \{(u_1,V)|~ |u_1|=1, u_1 \perp _{\C} V\}$. The spherical
subbundle $\tilde{H}^1 =\{ (u_1,u_2,V)| ~ |u_i|=1, u_1\perp _{\C}
u_2 \perp _{\C} V\}$ can be identified with $SU(n+2)/SU(n)$. By
using the diagonal action of the $S^1$ and $\tilde{\pi}$ one gets
the complex line bundle $\pi :H=\tilde{H}/S^1 \rightarrow
F=\tilde{F}/S^1$. The spherical subbundle in $\pi$ coincides with
the map $\pi^1: H^1=SU(n+2)/S[U(n)\times U(1)] \rightarrow F=
SU(n+2)/T$, where
$$
T= \left\{ \left( \begin{array}{ccc}  z & 0 & 0 \\ 0 & \bar{z}\det{\bar{A}} & 0 \\
0 & 0 & A \end{array} \right)| z\in U(1), A \in U(n) \right\}.
$$
Notice that $H^1$ is a $3$-Sasakian manifold and coincides with
Aloff-Wallach space for $n=1$ and $\pi ^1$ is its fibration over the
respective twistor space $\mathcal{Z}=F$ --- the space of complex
flags $\{([u],V)|~u \in S^{2n+3}, V \in G_n(\C^{n+2}), u \perp _{\C}
V \}$ in $\C^{n+2}$. It is not difficult to verify that the length
with respect to the metric $\bar{G}_{\alpha}$ of the characteristic
vector field dual to the form $\eta_1$ at the start time $r=1$ is
equal to $2(n+1)$. For metric $\bar{G}_{\alpha}$ to be well-defined
it is necessary for the circle generated by that vector field to be
factorized by the discrete subgroup $\Z _{n+1}$ because in the
$3$-Sasakian fibration $\pi^1$ there is already factorization by
$\Z_2$ (look \cite{Baz1} for details). Thus, the metrics
$\bar{G}_\alpha$ for $0\leq \alpha <1$ are defined on the tensor
power $\pi^{n+1}$. The theorem is proved.

\vskip8mm

\textbf{Acknowledgements.} I am grateful to Ya.V. Bazaikin for
discussions and to F. Aicardi for preparing this english version.

\end{document}